\documentstyle[12pt]{article}
\title{A Generalisation of Obata's theorem}
\author{Akhil Ranjan and G.~Santhanam}
\newtheorem{theorem}{Theorem}
\newtheorem{lemma}{Lemma}
\newtheorem{propn}{Proposition}

\newcommand{\ncmd}{\newcommand}
\ncmd{\minnu}{\mbox{$\parallel\!\nabla u\!\parallel$}}
\ncmd{\dt}{\mbox{${{\partial}\over{\partial t}}$}}
\ncmd{\dte}{\mbox{${{\partial}\over{\partial\theta}}$}}
\ncmd{\jt}{\mbox{${J_{\theta}}$}}
\ncmd{\euu}{\mbox{$E_{-{{u+1}\over 2}}$}}
\ncmd{\eu}{\mbox{$E_{-u}$}}
\ncmd{\orr}{\mbox{$\overline{R}$}}
\ncmd{\wm}{\mbox{$\widetilde{\sigma}(\theta)$}}
\ncmd{\ste}{\mbox{$\sigma(\theta)$}}
\ncmd{\nabu}{\mbox{$\nabla u$}}
\ncmd{\nabx}{\mbox{$\nabla _X$}}
\ncmd{\navu}{\mbox{$\nabla_{{\nabla u}}$}}
\ncmd{\oga}{\mbox{$\overline{\gamma}$}}
\ncmd{\olim}{\mbox{$\overline{m}$}}
\ncmd{\oliu}{\mbox{$\overline{u}$}}
\ncmd{\olimm}{\mbox{$\overline{M}$}}
\ncmd{\orc}{\mbox{$\overline{R}$}}
\ncmd{\ioga}{\mbox{$I_{\overline{\gamma}}$}}
\ncmd{\oho}{\mbox{$\overline{h}$}}
\ncmd{\grhu}{\mbox{$\nabla_hu$}}
\ncmd{\grhau}{\mbox{$\nabla_{\oho}u$}}
\begin{document}
\maketitle
\begin{abstract}
In a complete Riemannian manifold $(M, g)$ if the hessian of a real
valued function satisfies some suitable conditions then it restricts the
geometry of $(M, g)$. In this paper we characterize all compact rank-1
symmetric spaces, as those Riemannian manifolds $(M, g)$ admitting a
real valued function $u$ such that the hessian of $u$ has atmost two
eigenvalues $-u$ and $-{{u+1}\over 2}$, under some mild hypothesis on
$(M, g)$. This generalises a well known result of Obata which
characterizes all round spheres.
\end{abstract}
\footnote{1991 Mathematics Classification: 53C20, 53C22, 53C35}
\section{Introduction}
Lichnerowicz proved in \cite{L} that if $(M, g)$ is a
complete Riemannian manifold of dimension $n\geq 2$ such that the Ricci
tensor $Ric$ and the  metric $g$
verify the relation $Ric \geq l g$ for some $l >0$, then the first
eigenvalue $\lambda_1$ of the Laplacian of $(M, g)$ satisfies the
inequality
		$\lambda_1\geq {n\over {n-1}}l$.
While characterising the equality case of the above result,
Obata proved in \cite{MO} that a complete Riemannian manifold $(M, g)$ of
dimension $n\geq 2$ is isometric to the round sphere $(S^n, ds^2)$ if and
only if there is a real valued function $u\in C^2(M)$ such that, the
hessian of $u$, $\nabla^2u=-u Id$.

Recently Robert Molzon and Karen Pinney \cite{RK} have proved that
a complete K\"{a}hler manifold $(M, g, J)$ is isometric
to a complex projective space if and only if there is a real valued
function $u\in C^2(M)$ such that on $\{p\in M: \nabu(p)\neq 0\}$,
$\nabla^2u=-u Id
+({{u-1}\over 2})(Id-\Pi_{\grhu}-\Pi_{\grhau}$); here
$\grhu$ denotes the holomorphic gradient of $u$,
$\grhau$ denotes the antiholomorphic gradient of $u$ and
$\Pi_X$ denotes the orthogonal projection on the subbundle $X$. In their
paper the authors assume that $\nabu$ is an eigenvector of
$\nabla^2u$ with eigenvalue $-u$ and also that the mulitplicity of the
eigenvalue $-u$ is $2$. Since the manifold $(M, g, J)$ is assumed to be
Kahler, $J\nabu$ is also an eigenvector of $\nabla^2u$ with
eigenvalue $-u$ and hence the subbudle spanned by $\nabu$ and $J\nabu$
becomes a totally geodesic integrable subundle of $TM$.

In this paper, we drop the K\"{a}hler condition and also the condition
on the multiplicity of the eigenvalues of $\nabla^2u$. As a consequence,
we characterize all the compact rank-1 symmetric spaces under some mild
additional hypothesis. However in K\"{a}hler case we have much stronger
assertion. (See theorem 3).

Our method is different and uses both geometry and topology.
\begin{theorem}
Let $(M, g)$ be a complete Riemannian manifold of dimension $d$.
Let $u\in C^2(M)$ be a
real valued function such that the hessian of $u$, $\nabla^2u$, has atmost
two eigenvalues $-u$ and $-{{u+1}\over 2}$ and $\nabla u$ is an
eigenvector of $\nabla^2u$ with eigenvalue $-u$. Then
\begin{enumerate}
\item The multiplicity $k$ of the eigenvalue $-u$ is  $1$, $2$, $4$,
$8$ or $d$.
\item  If $k=1$, then either $(M, g)$ is isometric to
$I\!\!RI\!\!P^d$ or
$S^d$ with constant sectional curvature ${1\over 4}$.
\item If $k=$d, then $(M, g)$ is isometric to $S^d$ with
constant sectional curvature $1$. ({\bf Obata's theorem})
\item If $k=2$, $4$ or $8$, then $(M, g)$ is a pointed Blaschke manifold
at $m\in M$, where $m$ is the unique maxmium for the function $u$, with
totally geodesic cut locus $C(m)$. Moreover
$H^*(M, Z\!\!Z)=H^*(\olimm, Z\!\!Z)$ where $\olimm$ is a compact rank-1
symmetric space of dimension $kn$ and $k$ is the degree of the
generator of $H^(\olimm, Z\!\!Z)$.
\end{enumerate}
\end{theorem}
\begin{theorem}[Weak Obata's Theorem]
Let $(M, g)$, $u$, $m$ and $k\geq 2$ be as in theorem 1. Then
$Vol(M)=Vol(\olimm)$ where $\olimm$ is a compact rank-1 symmetric
space of dimension $kn$ with sectional curvature ${1\over 4}\leq
K_{\olimm}\leq 1$ and $H^*(M, Z\!\!Z)=H^*(\olimm, Z\!\!Z)$.
\end{theorem}
\begin{theorem}
Let $(M, g)$ and $u$ be as in theorem 1. If $(M, g)$ is a K\"{a}hler
manifold, then $(M, g)$ is isometric to $I\!\!\!CI\!\!P^n$ with sectional
curvature ${1\over 4}\leq K_{I\!\!\!CI\!\!P^n}\leq 1$
\end{theorem}
\begin{theorem}
Let $(M, g)$, $u$, $m$ and $k$ be as in theorem 1. If $k\geq 2$ and
$(M, g)$ is a $P_{2\pi}$ manifold at $m\in M$, then $(M, g)$ is isometric
to $\olimm$ where $\olimm$ is as in theorem 2.
\end{theorem}

\section{Preliminaries}
Let $(M, g)$ be a complete Riemannian manifold. For any $u\in C^2(M)$,
let $X:={\nabu\over{\parallel\!\nabu\!\parallel}}$ on
$\{ p\in M:\nabu(p)\neq 0\}$.
Then we have the following
\begin{propn}
If $\nabu$ is an eigenvector of $\nabla^2u$ then the integral curves of
$X$ are geodesics and conversely.
\end{propn}
\noindent{\bf Proof:}
\begin{eqnarray*}
\nabx X & = & {1\over{\parallel\!\nabu\!\parallel}}\nabx\nabu
						  +X({1\over{\parallel\!\nabu\!\parallel}})\nabu
\end{eqnarray*}
Then
\begin{eqnarray*}
\nabx X & = & 0
\end{eqnarray*}
iff
\begin{eqnarray*}
{1\over{\parallel\!\nabu\!\parallel}}\nabx\nabu & = &
						-X({1\over{\parallel\!\nabu\!\parallel}})\nabu \\
{} & = & {{X(\minnu)}\over{\minnu}^2}\nabu \\
{} & = & {{<\nabla_{\nabu}\nabu, \nabu>}\over{\minnu}^3}X \\
{} & = & {1\over{\parallel\!\nabu\!\parallel}}<\nabx\nabu, X> X
\end{eqnarray*}
iff
\begin{eqnarray*}
\nabla_{\nabu}\nabu & = & <\nabx\nabu, X> \nabu
\end{eqnarray*}
iff $\nabu$ is an eigenvector of $\nabla^2u$.
This completes the proof of the proposition.
\begin{propn}
Let $u\in C^2(M)$ be such that the integral curves of $X$ are geodesics.
Then $u$ does not have saddle points.
\end{propn}
\noindent{\bf Proof:} Suppose the proposition is false.

Let $p\in M$ be a saddle point of $u$. Then $\nabla^2u(p)$ has both
positive and negative eigenvalues. Hence there is a neighbourhood
$W$ of $p$
in $M$ such that the flow of $X$ have the form of hyperbolas near $p$ and
in this neighbourhood $W$ of $p$ they form a saddle. We may assume that
$W:=\exp_p(W_1)$ where $W_1$ is a
neighbourhood of $0$ in $T_pM$. (See \cite {A}, \cite{JM}).

Let $E^{us}\subseteq T_pM$ denote the eigensubspace
of $\nabla^2u(p)$ on which $\nabla^2u(p)$ is negative definite and
$E^s\subseteq T_pM$ denote the eigensubspace of $\nabla^2u(p)$ on which
$\nabla^2u(p)$ is positive definite. Let $W^{us}:=\exp_p(W_1\cap E^{us})$
and $W^s:=\exp_p(W_1\cap E^s)$. Then the integral curves of $X$ through any
point in $W^{us}$ will start from $p$ and diverge near $p$ in $W$ and
the integral curves of $X$ through any point in $W^s$ converge to $p$.

Let $B(p, \epsilon)$ be a strongly geodesically convex neighbourhood of
$p$ such that $B(p, \epsilon)\subseteq W$. Let
$T_{{\epsilon}\over 2^k}W^s$ be the tubular neighbourhood of radius
${{\epsilon\over 2^k}}$ of $W^s$ and $T_{{\epsilon}\over 2^k}W^{us}$ be the
tubular neighbourhood of radius ${{\epsilon}\over 2^k}$ of $W^{us}$ for
$k\geq 1$.
Now we choose a point
$q_1\in B(p, \epsilon)\cap(T_{{\epsilon}\over 2^k}W^s\setminus W^s)$. Let
$\sigma$ be the minimizing geodesic from $q_1$ to
$B(p, \epsilon)\cap W^s$ such that $\sigma(0)=q_1$ and
$\sigma(1)=q_2\in B(p, \epsilon)\cap W^s$.
Let $\gamma_s$ denote the integral curve of $X$, starting at $\sigma(s)$.
If $k$ is large, these geodesics $\{\gamma_s\}$ will pass through the
tubular neighbourhood $T_{{\epsilon}\over 2^k}W^{us}$ of $W^{us}$
and these geodesics will converge to a geodesic in piecewise $C^1$ limit.
The limiting
geodesic will pass through $p$ and broken at $p$.
Since the geodesics $\gamma_s$ are
all minimizing in $B(p, \epsilon)$, the limiting geodesic
will also be a minimizing geodesic in
$B(p, \epsilon)$. This is a contradiction. Hence $u$ can't have saddle
points.  This completes the proof of the proposition.

Since $u$ does not have saddle points, the only possible critical points
of $u$ are maxima and minima. To describe these points we first compute
the function $u$ along the integral curves of $X$ in the following
\begin{lemma}
Let $u\in C^2(M)$ be such that $\nabu$ is an eigenvector of $\nabla^2u$
with eigenvalue $-u$. Then
\begin{enumerate}
\item Along the integral curves $\gamma$ of $X$, the function $u$ is of
the form $u(\gamma(t))=A_{\gamma}\cos t + B_{\gamma}\sin t$.
\item The functiuon $u$ has only isolated critical points along $\gamma$.
\end{enumerate}
\end{lemma}
\noindent{\bf Proof:} Let $\gamma$ be an integral curve of $X$.
Since $\nabu$ is an eigenvector with eigenvalue $-u$, we have that
\begin{eqnarray*}
u''(\gamma(t)) & = & <\nabla_{\gamma'(t)}\nabu, \gamma'(t)> \\
{} & = & -u(\gamma(t)) <\gamma'(t), \gamma'(t)>
\end{eqnarray*}
whenever $\nabu(\gamma(t))\neq 0$.
Hence $u(\gamma(t))= A_{\gamma}\cos t+ B_{\gamma}\sin t$ whenever
$\nabu(\gamma(t))\neq 0$. Since $(M, g)$ is a complete Riemannian manifold,
any geodesic can be
extended for all $t\in I\!\!R$ and we can write
$u(\gamma(t))= A_{\gamma}\cos t+ B_{\gamma}\sin t$
whenever $\nabu(\gamma(t))\neq 0$.

Clearly critical points of $u$ are not isolated along $\gamma$ only if
$A_{\gamma}=B_{\gamma}=0$. Since $\gamma'(t)=X(\gamma(t))$ for almost all
$t$, this will mean that $X=0$, a contradiction for a unit vector field.
This proves the lemma.

Here afterwards we will write $u(t)$ for $u(\gamma(t))$ and $X(t)$ for
$X(\gamma(t))$.
\begin{lemma}
Let $u$ be as in lemma 1. Then
\begin{enumerate}
\item The function $u$ attains its maximum at some point $m\in M$ and
$\nabla^2u$ is non-degenerate at $m$.
\item $u(q)=\cos d(q, m)$ for any point $q\in M$.
\end{enumerate}
\end{lemma}
\noindent{\bf Proof:} Let $\gamma$ be an integral curve of $X$. We know
from lemma 1 that $u(t)= A_{\gamma} \cos t + B_{\gamma} \sin t$ for all
$t\in I\!\!R$.
Therefore the function $u$ attains a positive maximum and a negative
minimum along the geodesic $\gamma$. We may assume that the function $u$
attains its maximum along $\gamma$ at $t=0$. Let $\gamma(0)=m$. Clearly
$u(m)>0$. Since $u'(0)=0$, we have that $B_{\gamma}=0$.

Since $\gamma$ is an integral curve of $X$,
		$\gamma'(t) = X(t) =  {{\nabu(t)}\over
									{\parallel\nabu(t)\parallel}}$
whenever
			$\nabu(t)\neq 0$.
Therefore
\begin{eqnarray*}
u'(t)  & = & <\nabu(t), \gamma'(t)> \\
{}     & = &  \parallel\nabu(t)\parallel
\end{eqnarray*}
whenever $\nabu(t)\neq 0$. Hence $u'(t)=0$ iff
$\parallel\nabu(t)\parallel=0$.
This shows that $t=0$ is a critical point for the function $u$.
Since $u(m)>0$ and the hessian of $u$ has atmost two eigenvalues
$-u$ and $-{{u+1}\over 2}$, we have that $\nabla^2u(m)$ is negative
definite at $m$. Hence $\nabla^2u$ is non-degenerate at $m$.

Since $\nabla^2u(m)$ is negative definite,
$m$ is an isolated critical point. Therefore there is a neighbourhood $W$
of $m$ such that the integral curves of $X$ passing through the points in
$W$ will start at $m$. This proves that any geodesic $\gamma$ starting at
$m$ is tangent to $\nabu$ and hence any $v\in T_mM$ is an
eigenvector of $\nabla^2u(m)$ with eigenvalue $-u(m)$.
Therefore all the eigenvalue of
$\nabla^2u(m)$ must be equal and hence $A_{\gamma}=1$ for all geodesics
$\gamma$ starting at $m$ and $u(t)=\cos t$ along any  geodesic $\gamma$
starting at $m$.

Since $(M, g)$ is a complete Riemannian manifold, for any point $q\in M$
there is a distance minimizing geodesic $\gamma$ from $m$ to $q$.  Hence
$u(q)=\cos d(m, q)$. This completes the proof of the lemma.

Now we describe the set $C:=\{q\in M: u(q)=\min_{p\in M}u(p)\}$ in the
following
\begin{lemma}
The function $u$ has atmost two maximum $m$ and $m'$.
\begin{enumerate}
\item If $m\neq m'$, then the cut locus of $m=\{m'\}$ and the cut locus
of $m'=\{m\}$. Moreover $M$ is homeomorphic to $S^{{\rm dim} M}$ and the
multiplicity of the eigenvalue $-u$ is $1$.
\item
\begin{enumerate}
\item If $m'=m$, then $(M, g)$ is a pointed Blaschke manifold at
$m\in M$, $C$ is the cut locus of $m$ and $C:=\{q\in M : d(q, m)=\pi\}$.
Further $C$ is totally geodesic.
\item The multiplicity $k$ of the eigenvalue $-u$ is $1$, $2$, $4$, $8$ or
dim$M$.
\item If $k=1$ then $\pi_1(M)=Z\!\!Z_2$ and $M$ has the homotopy type of
$I\!\!RI\!\!P^d$. If $k= d$, then $M$ has the
homotopy type of $S^d$. If $k=2$, $4$ or $8$ then $M$ is simply
connected and $H^*(M, Z\!\!Z)= H^*(\olimm, Z\!\!Z)$ where $\olimm$ is as
in theorem 1.
\end{enumerate}
\end{enumerate}
\end{lemma}
\noindent{\bf Proof:} Let $C_M:=\exp_m(S(0, 2\pi))$, the image of the
sphere of radius $2\pi$ in $T_mM$. Then any point $q\in C_M$ is a point of
maximum for the function $u$. i.e., $u(q)=1$ for all $q\in C_M$. Since
$C_M$ is connected and $\nabla^2u(q)$ is non-degenrate for $q\in C_M$, we
have that $C_M$ is a singleton. Let $C_M=\{m'\}$. This proves that the
function $u$ has atmost two points of maximum and also that $-1$ is the
minimum for the function $u$. Hence $C:=\{q\in M : d(q, m)=\pi\}$
and the eigenvalues of $\nabla^2u(q)$ are atmost $-1$ and
$-{{u(q)+1}\over 2}=0$ for all $q\in C$.

We have seen above that $C=\exp_m(S(0, \pi)$. Since $C$ is connected and
$\nabla^2u$ has atmost two eigenvalues $-1$ and $0$ on $C$, the rank of
$\nabla^2u$ is constant on $C$. Let us denote this constant by $k$. Then
$C$ is a $(d-k)$- dimensional submanifold of $M$ and the normal bundle of
$C$ is spanned by limiting vectors ${{\nabu}\over{\minnu}}$ as we move
towards $C$.

Since $\nabu\neq 0$ on the open geodesic ball $B(m, \pi)\setminus \{m\}$
and the flow of $\nabu$ are geodesics, for any geodesic $\gamma$ starting
at $m$, the cut point of $m$ along the geodesic $\gamma$ does not occur on
$B(m, \pi)\setminus \{m\}$. Hence $\exp_m\mid_{B(0, \pi)} : B(0, \pi)\to M$
is a diffeomorphism of the open ball $B(0, \pi)$ of radius $\pi$ in $T_mM$
on to the geodesic ball $B(m, \pi)$.

If $m'\neq m$, then, since $\nabu\neq 0$ on $B(m', \pi)\setminus\{m'\}$
and the flow of $\nabu$ are geodesics, the set $B(m', \pi)\setminus\{m'\}$
will be free of cut points of $m$.  Therefore the cut locus of $m$ is
contained in $C\cup\{m'\}$.
Now, since $M=\overline{B(m, \pi)}\cup_C\overline{B(m', \pi)}$ and $M$ is
a smooth manifold $C$ is an $(d-1)-$ dimensional submanifold of $M$.
Therefore $\exp_m : S(0,\pi)\to C$ is either one-one or a two sheeted
covering. If $\exp_m : S(0, \pi)\to C$ is a two sheeted covering, then all
the geodesics starting at $m$ will stop minimizing beyond $C$. This is a
contradiction.  Hence $\exp_m: S(0, \pi)\to C$ is one-one. This shows that
diam$M>\pi$.  Now the flow of $\nabu$ will move towards $m$ for points $q$
at distance $<\pi$ from $m$ and the flow of $\nabu$ will move towards
another maximum $m'$ for points $q$ at distance $>\pi$ from $m$. This
proves that $M$ is homeomorphic to $S^d$. The proof also shows
that the cut-locus of $m$ is $\{m'\}$ and the cut-locus of $m'$ is
$\{m\}$.

Since $C$ is a submanifold of dimension $(d-1)$, the multiplicity $k$ of
the eigenvalue $-u$ is $1$.

If $m'=m$, then $C$ is the cut locus of $m$.i.e., $C=C(m)$.
Since the tangential cut
locus of $m$ is spherical, it follows from \cite{O}, \cite{NS} that
$(M, g)$ is a Blaschke manifold at $m\in M$. (See also \cite{BE}).
Now we prove that $C$ is totally geodesic.

Let $v\in T_qC$. We extend $v$ to a vectorfield $V$ in a neighbourhood of
$q\in M$ such that $V$ is tangential to the level sets of $u$. Now we
write $V=V_1+V_2$ where $V_1$ is an eigenvectorfield of $\nabla^2u$ with
eigenvalue $-u$ and $V_2$ is an eigenvectorfield of $\nabla^2u$ with
eigenvalue $-{{u+1}\over 2}$. Then
\begin{eqnarray*}
\nabla_VX & = & \nabla_{V_{1}}X+\nabla_{V_{2}}X \\
{}        & = & -{{u}\over\minnu}V_1-{{u+1}\over{2\minnu}}V_2
\end{eqnarray*}
Therefore, since $V_2(q)=v$ and $V_1(q)=0$, we have that
\begin{eqnarray*}
<\nabla_VX, V>(q) & = &  -{{u(q)}\over\minnu}{\parallel V_1(q)\parallel}^2
								-{{u(q)+1}\over{2\minnu}}{\parallel V_2(q)\parallel}^2\\
{}                & = & 0
\end{eqnarray*}
This proves that $C$ is a totally geodesic submanifold of $(M, g)$.

Since $(M, g)$ is Blaschke manifold at $m$, it follows from
Bott-Samelson's theorem \cite{BE}, that all geodesics starting at $m$ have
same index $\lambda=0$, $1$, $3$, $7$ or $d-1$.
If $\lambda>0$, then we have the following possiblities.
\begin{enumerate}
\item $\lambda=1$, $d=2n$ and $M$ has the homotopy type of
$I\!\!\!CI\!\!P^n$.
\item $\lambda=3$, $d=4n$ and $M$ has the integral cohomolgy ring of
$I\!\!HI\!\!P^n$.
\item $\lambda=7$, $d=16$ and $M$ has the integral cohomology ring of
$I\!\!\!CaI\!\!P^2$.
\item $\lambda= d-1$ and $M$ has the homotopy type of $S^d$.
\end{enumerate}
If $\lambda=0$, then $\pi_1(M)=Z\!\!Z_2$ and $M$ has the homotopy type of
$I\!\!RI\!\!P^d$.

When $\lambda>0$, the cut locus of $m$ coincides with the conjugate locus
of $m$.
Since the cut-locus coincides with the conjugate locus, we have that
$\lambda=k-1$, where $k$ is the rank of the hessian of
$u$, $\nabla^2u$ on $C$. This proves that the multiplicity of the
eigenvalue $-u$ is $1$, $2$, $4$, $8$ or dim$M$. This completes the proof.

\noindent{\bf Remark:} If $k = d$, then $-{{u+1}\over 2}$ is not an
eigenvalue of $\nabla^2u$. Hence $\nabla^2u=-u Id$ and $C$ is singleton.
This is Obata's theorem.

\section{Proof of theorem 1}
\noindent{\bf Proof of 1(1) and 1(4):}
Proof of theorem 1(1) and theorem 1(4) follows from lemma 3.

\noindent{\bf Proof of 1(2):}
Now we prove that if $k=1$, then either $(M, g)$ is isometric to $S^d$ or
$I\!\!RI\!\!P^d$ with constant sectional curvature ${1\over 4}$.

Since the multiplicity of the eigenvalue $-u$ is $1$, any vector
$E\perp\nabu$ is an eigenvector of $\nabla^2u$ with eigenvalue
$-{{u+1}\over 2}$. Hence the eigensubbundle
$\euu:=\{E\in TM : \nabla^2u(E)=-{{u+1}\over 2}E\}$ is parallel along $X$.

Let $\gamma$ be a geodesic starting at $m$ and let $J$ be the Jacobi
field
describing the variation of the geodesic $\gamma$ such that $J(0)=0$ and
$J'(0)=E\in\euu$ is of unit norm. Since $[J, X]=0$ along
$\gamma$, we have that $\nabla_XJ=\nabla_JX$.
Since $u=\cos t$ along the geodesics $\gamma$ starting at $m$,
$\nabu=-\sin t {{\partial}\over{\partial t}}$ and
$X={{\nabu}\over{\minnu}}=-{{\partial}\over{\partial t}}$.
Therefore
\begin{eqnarray*}
-J' & = & {1\over \minnu}\nabla_J\nabu \\
{}  & = & -{{u+1}\over 2}{1\over\minnu}J \\
{}  & = & -{1\over 2}{{\cos{t\over 2}}\over{\sin{t\over 2}}} J
\end{eqnarray*}
and
\begin{eqnarray*}
{{<J', J>}\over{{\parallel J\parallel}^2}} & = &
		{1\over 2}{{\cos{t\over 2}}\over{\sin{t\over 2}}}
\end{eqnarray*}
Therefore
\begin{eqnarray*}
{d\over dt}\log{{\parallel J\parallel}\over{\sin{t\over 2}}} & = & 0
\end{eqnarray*}
and hence
${{\parallel J\parallel}\over{\sin{t\over 2}}} =
		{{\parallel J\parallel}\over{\sin{t\over 2}}}\mid_{t=0}=2$.
This shows that
		$\parallel J\parallel =2\sin{t\over 2}$.
Since $\euu$ is parallel along
$X$, we can write
		$J(t)=2\sin{t\over 2}E(t)$
where $E$ is a parallel vectorfield along $X$.
Therefore
\begin{eqnarray*}
R(J, X)X & = & -J'' \\
{}       & = & {1\over 4}J
\end{eqnarray*}
This proves that $\euu$ is an eigensubundle of $R(\cdot, X)X$ with
eigenvalue ${1\over 4}$ and
$K(E, X)={1\over 4}$ for $E\in\euu$ of unit norm.

First we prove that if $M$ is homeomorphic to $S^d$, then $(M, g)$ is
isometric to $S^d$ with constant sectional curvature ${1\over 4}$.

We choose a point $\olim\in S^d$ and fix a linear isometry
$i : T_mM\to T_{\olim}S^d$. Now we define a map $\Phi: M\to S^d$ by
$\Phi(q):=\exp_{\olim}\circ i\circ\exp_m^{-1}(q)$. Then $\Phi$
maps the geodesics $\gamma$ starting at $m$ in $M$ on to the geodesics
$\overline{\gamma}$ starting at $\olim$ in $S^d$ and it also maps the
geodesic spheres around $m$ in $M$ on to the geodesic spheres around
$\olim$ in $S^d$. To complete the proof, we only have to show that
$d\Phi$ is norm preserving. But this follows easily from the fact that
the Jacobi field along the geodesics $\gamma$ atarting at $m$ in $M$ are
same as that of the Jacobi fields along the geodesics $\overline{\gamma}$
starting $\olim$ in $S^d$. This completes the proof.

When $M$ has the homotopy type of $I\!\!RI\!\!P^d$, a similar argument as
above shows that $(M, g)$ is isometric to $I\!\!RI\!\!P^d$ with constant
sectional curvature ${1\over 4}$.

\noindent{\bf Proof of 1(3):} Since $k=d$, we have that $\nabla^2u=-uId$.
Now an argument similar to the proof of 1(3) shows that $M$ is isometric
to $S^d$ with constant sectional curvature $1$.
\section{Proof of theorem 2}
Let $S(m, r):=\{q\in M: u(q)=\cos r\}$. Then
$L:={{\nabla^2u}\over{\minnu}}$ is the second fundamental form the level
sets $S(m, r)$ of the function $u$, with respect to the inward unit
normal.(See \cite{HK}).
Hence the mean curvature of $S(m, r)$ at any point $p\in S(m, r)$ is
\begin{eqnarray*}
Tr(L(p)) & = & \sum_{i=1}^{kn-1}<L(p)(e_i), e_i>
\end{eqnarray*}
where $e_1$, $e_2, \cdots , e_{kn-1}$ is an orthonormal basis of
$T_pS(m, r)$. Since $-u$ and $-{{u+1}\over 2}$ are the only eigenvalues of
$\nabla^2u$ and the multiplicity of the eigenvalue $-u$ is $k$, we have
that
\begin{eqnarray*}
Tr(L(p)) & = & \sum_{i=1}^{kn-1}<L(p)(e_i), e_i> \\
{} & = & -\big\{(k-1)\cot t +
					({{kn-k}\over{\sin t}}){{1+\cos t}\over 2}\big\} \\
{} & = & -\big\{(k-1)\cot t + {{kn-k}\over 2}\cot{t\over 2}\big\}
\end{eqnarray*}
On the other hand, we know that
$Tr(L(p))={{\theta_m'(t)}\over{\theta_m(t)}}$ where $\theta_m(t)$ is the
Riemannian volume density function of $(M, g)$ in geodesic polar
coordinates centred at the point $m$ in $M$.(See \cite{HK}). This proves
that
\begin{eqnarray*}
\theta_m(t) & = & 2^{kn-k}\sin^{kn-k}{t\over 2}\sin^{k-1}t
\end{eqnarray*}
Therefore
\begin{eqnarray*}
Vol(M) & = & \int_{U_{m}M}\int_0^{\pi}\theta_m(t)dt\, d\theta \\
{}     & = & \int_{U_{m}M}\int_0^{\pi}
						 2^{kn-k}\sin^{kn-k}{t\over 2}\sin^{k-1}t \\
{}     & = & Vol(\olimm)
\end{eqnarray*}
This completes the proof of theorem 2.
\section{Proof of Theorem 3}
First we note that, since $(M, g)$ is a K\"{a}hler manifold, the even
betti numbers are positive.(See \cite{GH}). Hence $H^2(M, Z\!\!Z)\neq 0$.
This proves that the multiplicity of the eigenvalue $-u$ is
$2$ and $H^2(M, Z\!\!Z)=H^2(I\!\!\!CI\!\!P^n, Z\!\!Z)$.

For each point $q\in C$, the cut-locus of $m$, we denote by $\sum_q$, the
union of all geodesics from $m$ to $q$. Then $\sum_q$ is a smooth surface
except possibly at $m$ and $\sum_q$ is totally geodesic at $q$.
(See \cite{NS}). We write the induced metric on $\sum_q$ by
\begin{eqnarray*}
ds^2 & = & dr^2+f_q(r, \theta)d\theta^2
\end{eqnarray*}
where $f_q(r, \theta)$ is a continuous function which is smooth except
possibly at $m$ and $f_q(0, \theta)=0=f_1(\pi, \theta)$ and
$f_q'(\pi, \theta)=-1$.
We, now, prove the following
\begin{lemma}
\begin{enumerate}
\item $\sum_q$ is a smooth totally geodesic surface in $M$.
\item $\sum_q$ is isometric to $S^2$ with constant curvature $1$.
\end{enumerate}
\end{lemma}
\noindent{\bf Proof:} Let $\gamma$ be a geodesic segment joining $m$ and
$q$ and $J$ be the Jacobi field describing the variation of $\gamma$ such
that $J(0)=0$ and $J(\pi)=0$. We normailse $J$ such that
$\parallel\!J'(\pi)\!\parallel=1$. Since $J(t)\subseteq T\sum_q$ and
$\sum_q$ is totally geodesic at $q$, $J'(\pi)\in T\sum_q$.

Since $J$ is a Jacobi field along $\gamma$, we have that $[J, X]=0$ along
$\gamma$. Hence
\begin{eqnarray*}
\nabla_XJ & = & \nabla_JX
\end{eqnarray*}
and
\begin{eqnarray*}
-<J', J> & = & {1\over {\minnu}}<\nabla_J\nabu, J>
\end{eqnarray*}
Since the eigenvalues $-{{u+1}\over 2}$ and $-u$ of $\nabla^2u$ satisfy
the inequality $-{{u+1}\over 2}\leq -u$, it follows that
\begin{eqnarray*}
-<J', J> & \leq & -{u\over{\minnu}}\parallel\!J\!\parallel^2
\end{eqnarray*}
and
\begin{eqnarray*}
{{<J', J>}\over {\parallel\!J\!\parallel^2}} & \geq &
							{{\cos t}\over {\sin t}}
\end{eqnarray*}
Therefore
\begin{eqnarray*}
{d\over{dt}}\log{{\parallel\! J\!\parallel}\over{\sin t}} & \geq & 0
\end{eqnarray*}
This shows that
\begin{eqnarray*}
{{\parallel\! J\!\parallel}\over{\sin t}} & \leq &
		{{\parallel\! J\!\parallel}\over{\sin t}}\!\mid_{t=\pi} \\
{} & = & 1
\end{eqnarray*}
i.e., ${{\parallel J\parallel}\over{\sin t}}\leq \sin t$. Hence
\begin{eqnarray*}
Area(\sum_{q}) & = & \int_{S_{m}\sum_{q}}\int_0^{\pi}
		{\parallel\! J(t)\!\parallel} dt\, d\theta \\
{} & \leq & \int_{S_{m}\sum_{q}}\int_0^{\pi} \sin t \,dt\, d\theta \\
{} & = & 4\pi
\end{eqnarray*}
This proves that $Area(\sum_q)\leq 4\pi$. Therefore
\begin{eqnarray*}
\min_{q\in C}Area(\sum_q) & \leq & 4\pi
\end{eqnarray*}
and equality holds iff
\begin{enumerate}
\item The Jacobi field $J$ is an eigenvector field of $\nabla^2u$ with
eigenvalue $-u$, $\parallel\!J\!\parallel =\sin t$ and $\sum_q$ is
a smooth totally geodesic surface in $M$.
\item $\sum_q$ is isometric to $S^2$ with sectional curvature $1$.
\item $\exp_m : S(0, \pi)\to C$ is a great circle fibration; here
$S(0, \pi)$ denotes the sphere of radius $\pi$ in $T_mM$.
\end{enumerate}
We will now show, using K\"{a}hler condition, that $Area(\sum_q)=4\pi$
for each $q\in C$.

Let $\omega$ be the K\"{a}hler form of $(M, g)$. Then we know that
${{\omega}^n\over {n!}}$ is the volume form of $(M, g)$.(See \cite{BGM}).
On the other hand we know that
$Vol(M)=Vol(\olimm)={{4\pi^n}\over{n!}}$. Therefore
\begin{eqnarray*}
\int_m{{\omega^n}\over{n!}} & = {{4\pi^n}\over {n!}}
\end{eqnarray*}
i.e., $\int_M({{\omega}\over{4\pi}})^n=1$.

This shows that $({{\omega}\over{4\pi}})^n$ is a generator of
$H^{2n}(M, Z\!\!Z)$. Now we write ${{\omega}\over{4\pi}}=c\theta$ where
$c>0$ and $\theta$ is a generator of $H^2(M, Z\!\!Z)$. Since
$\int_M({{\omega}\over{4\pi}})^n=1$ and $\theta^n$ is a generator of
$H^{2n}(M, Z\!\!Z)$, we get that $c=1$. Thus we have proved that
${{\omega}\over{4\pi}}$ is a generator of $H^2(M, Z\!\!Z)$. Hence
\begin{eqnarray*}
\int_{\sum_{q}}{{\omega}\over{4\pi}} & = & 1
\end{eqnarray*}
But from Wirtinger's inequality (See \cite{FRH}), it follows that
$$
{{Area(\sum_q)}\over{4\pi}}\geq \int_{\sum_{q}}({{\omega}\over{4\pi}})=1
$$
This proves that $Area(\sum_q)=4\pi$. Hence the proof of the lemma.

Let $E_{-u}:=\{E\in TM: \nabla^2u(E)=-uE\}$. Then from lemma 4, we see
that $E_{-u}$ is parallel along $\gamma$ for any geodesic $\gamma$
from $m$ to $q$ and $E_{-u}\mid\!_{\sum_{q}}$ is the tangent bundle of
the surface $\sum_q$ for any $q\in C$.

Let $J$ be the Jacobi field along a geodsic $\gamma$ from $m$ to $q$ such
that $J(0)=0$ and $J(\pi)=0$. Then, since $\sum_q$ is isometric to $S^2$
with constant curvature $1$, the Jacobi field $J$ is of the
form $J(t)=\sin t E(t)$ where $E(t)$ is parallel
along $\gamma$ and $E(t)\in E_{-u}$ . This shows that
\begin{eqnarray*}
R(J, X)X & = & -J'' \\
{} & = & J
\end{eqnarray*}
This proves that $E_{-u}$ is an eigensubbundle of $R(., X)X$ with
eigenvalue $1$.

Since $E_{-u}$ is parallel along $X$, the eigensubbundle
$E_{-{{u+1}\over 2}}:=\{E\in TM:\nabla^2u(E)=-{{u+1}\over 2}\}$ is also
paralle along $X$. An easy computation shows that $E_{-{{u+1}\over 2}}$
is also eigensubbundle of $R(., X)X$ with eigenvalue $1\over 4$.
This shows that any Jacobi field $J$ describing the variation of a geodesic
$\gamma$ such that $J(0)=0$ and $J'(0)\in\euu$ is given by
$J(t)=2\sin{t\over 2}E(t)$ where $E(t)$ is a parallel field along $X$.
(See also proof of theorem 1 in section 3).

Now we prove the following
\begin{lemma}
$\exp_m: S(0, \pi)\to C$ is congruent to Hopf fibration.
\end{lemma}
{\bf Proof:} We know from lemma 4 that $\exp_m : S(0, \pi)\to C$ is a
great circle fibration. Here we will show that this fibration is Riemannian.
Then it will follow from the classification of Riemannian submersions of
round spheres with connected totally geodesic fibres
that this fibration is congruent to
Hopf fibration. (See \cite{ESC}, \cite{R}, \cite{GWZ}).

Let $W\in T_qC$ be a unit vector. Let $\gamma$ be a geodesic from $m$ to
$q$. Then, since $\euu$ is parallel along $\gamma$, there exists a
parallel vectorfield $E(t)$ along $\gamma$ such that $E(0)=E$ and
$E(\pi)=W$. Now the unit vector $E$ is tangential to $S(0, \pi)$ and
orthogonal to the fibre through $\pi\gamma'(0)=\pi v $ where
$\gamma'(0)=v$.

Let $J(t):=d(\exp_m)_{tv}(tE)$. Then $J$ is a Jacobi field along $\gamma$
such that $J(0)=0$ and $J'(0)=E$. But we have seen above that any such
Jacobi field is given by $J(t)=2\sin{t\over 2}E(t)$. Therefore
$d(\exp_m)_{\pi v}(\pi E)= J(\pi)=W$. This proves that, up to a constant
factor $\pi$, $\exp_m : S(0, \pi)\to C$ is a Riemannian submersion with
totally geodesic fibres and hence $\exp_m : S(0, \pi)\to C$ is congruent
to Hopf fibration.

Now we come to the proof of theorem 3.

\noindent{\bf Proof of theorem 3:}
Let us fix a point $\olim\in I\!\!\!CI\!\!P^n$. Then we know that
$\exp_{\olim} : S(0, \pi)\to C(\olim)$ is the standard Hopf fibration;
here $C(\olim)$ denotes the cut locus of $\olim$. Now
since,
$\exp_m : S(0, \pi)\to C$ is congruent to Hopf fibration there is a
linear isometry
		$i : T_mM\to T_{\olim}I\!\!\!CI\!\!P^n$
such that $i$ carries the fibres of the fibration
$\exp_m : S(0, \pi)\to C$ to the fibres of
$\exp_{\olim} : S(0, \pi)\to C(\olim)$.

Now we define
		$\Phi: M\to I\!\!\!CI\!\!P^n$
by
		$ \Phi(q) :=\exp_{\olim}\circ i\circ \exp_m^{-1}(q)$.
Then for any geodesic $\gamma$ through $m\in M$, $\oga:=\Phi(\gamma)$ is a
geodesic through $\olim\in I\!\!\!CI\!\!P^n$. To complete the proof, we
have to show that $d\Phi$ preserves the lengths of Jacobi fields along
$\gamma$.

Let $\gamma$ be a geodesic starting at $m\in M$. Let $\gamma'(0)=v$ and
let $E(t)$ be the parallel vectorfield along $\gamma$ such that
$E(t)\in\eu$. Then $E(t)$ is given by
$E(t)={1\over{\sin t}}d(\exp_m)_{tv}(E(0))$ and $d\Phi_{\gamma(t)}$ maps
$d(\exp_m)_{tV}(E(0))$ to $d(\exp_m)_{ti(V)}(i(E(0)))$. Since the
isometry preserves the fibration, the vector
$d(\exp_m)_{ti(V)}(i(E(0)))\in
		E_1:=\{w\in TI\!\!\!CI\!\!P^n: \orr(w, \oga')\oga'=w\}$; here
$\orr$ denotes the Riemannian curvature tensor of $I\!\!\!CI\!\!P^n$.
We know all the Jacobi fields on $I\!\!\!CI\!\!P^n$ along the geodesics
$\oga$ and they are of the form
$J(t)=\sin t E(t)$ for $E\in E_1$ and $J(t)=2\sin{t\over 2}E(t)$ for
$E\in E_{{1\over 4}}$ where
$E_{{1\over 4}}:=
\{ w\in TI\!\!KI\!\!P^n : \orr(w, \oga')\oga'={w\over 4}\}$. Hence we see
that $d(\exp_m)_{ti(V)}(i(E(0)))={{\sin t}\over t}i(E(0))$. This shows
that $d\Phi$ is norm preserving on $\eu$.

By similar arguments we can show that $d\Phi$ is norm preserving on
$\euu$. Hence $\Phi : M\setminus C\to I\!\!\!CI\!\!P^n\setminus C(\olim)$
is an isometry. Now it follows from the uniform continuity that
$\Phi : M\to I\!\!\!CI\!\!P^n$ is an isometry. This completes the proof of
theorem 3.
\section{Proof of theorem 4}
In this section we assume that $(M, g)$ is a $P_{2\pi}$ manifold at $m$
and $k\geq 2$.

We define the eigensubbundle $\eu:=\{E\in TM : \nabla^2u(E)=-uE\}$. First
we prove that $\eu$ and $\euu$ are parallel along $X$ and from this we
deduce that $\eu$ and $\euu$ are also eigensubbundles of $R(\cdot, X)X$ with
eigenvalues $1$ and ${1\over 4}$ respectively in the following
\begin{lemma}
\begin{enumerate}
\item $\exp_m : S(0, \pi)\to C$ is a great sphere fibration.
\item $\eu$ and $\euu$ are parallel along $X$.
\item $R(E, X)X=E$ if $E\in\eu$ and $R(E, X)X={E\over 4}$ if $E\in\euu$.
\end{enumerate}
\end{lemma}
{\bf Proof:} Since $k\geq 2$, we know from lemma 3 that the cut locus $C$
of $m$ coincides with the conjugate locus of $m$ and the multiplicity of
each conjugate point $q\in C$ is $k-1$.

Let $\gamma$ be a geodesic starting at $m$ and let $J$ be a Jacobi field
describing the variation of the geodesic $\gamma$ such that
$J(0)=0=J(\pi)$ and $\parallel J'(\pi)\parallel =1$.
Let $\gamma(\pi)=q\in C$ and ${\rm seg}(q, m)$ denote the set of all
minimizing geodesics from $q$ to $m$. Since $(M, g)$ is a Blaschke
manifold at $m$ and $C$ is the cut locus of $m$, we know from
Omori\cite{O}, Nakagawa-Shiohama\cite{NS} that
$\Lambda(m, q):=\{\gamma'(0) : \gamma\in {\rm seg}(q, m)\}$, the link of
$q$ and $m$, is a great subsphere of $T_qM$ orthogonal to $T_qC$. We denote
by $\sum_q$, the union of all geodesics from $q$ to $m$. Then $\sum_q$ is
a smooth $k$-dimensional submanifold of $M$ except possibly at $m\in M$
and $\sum_q$ is totally geodesic at $q$. Since $J(t)\subseteq T\sum_q$ and
$\sum_q$ is totally geodesic at $q$, we have that $J'(\pi)\in T_q\sum_q$.

Let $\ste:=\cos\theta \gamma'(\pi) +\sin\theta J'(\pi)$, the great
circle in the plane spanned by $\gamma'(\pi)$ and $J'(\pi)$. Let
$\gamma_{\theta}$ denote the geodesic from $q$ to $m$ such that
$\gamma_{\theta}'(0)=\ste$. Since $(M, g)$ is a $P_{2\pi}$ manifold at
$m\in M$, these geodesics $\gamma_{\theta}$'s are all smoothly closed at
$m$. Let $\wm$ be a curve in $U_mM$ defined by
$\wm:=\gamma_{\theta}'(\pi)$.

Now we have a variation $H(t, \theta):= \exp_m(t\wm)$ of the geodesics
$\gamma_{\theta}(t):=\exp_m(t\wm)$ for a fixed $\theta$. Then the Jacobi
field $J_{\theta}$ along $\gamma_{\theta}$ is given by
$J_{\theta}(t)=
		{{\partial}\over{\partial\theta}}\mid_{(t, \theta)}H(t, \theta)$.
Since $[J_{\theta}, X]=0$ along each $\gamma_{\theta}$, we have that
\begin{eqnarray*}
\nabla_XJ_{\theta} & = & \nabla_{J_{\theta}}X
\end{eqnarray*}
and
\begin{eqnarray*}
-<J_{\theta}', J_{\theta}> & = & {1\over\minnu}
																 <\nabla_{J_{\theta}}\nabu, J_{\theta}>
\end{eqnarray*}
Since the eigenvalues $-u$ and $-{{u+1}\over 2}$ of $\nabla^2u$ are such
that $-{{u+1}\over 2}\leq -u$, it follows, as in lemma 4, that
\begin{eqnarray*}
{{\parallel\jt\parallel}\over{\sin t}}\mid_{t=0} & \leq &
			{{\parallel\jt\parallel}\over{\sin t}}\mid_{t=\pi} \\
{} & = & 1
\end{eqnarray*}
i.e., $\parallel J_{\theta}'(0)\parallel\leq 1$.
Now
\begin{eqnarray*}
J_{\theta}'(0) & = & \dt\mid_{t=0}\dte\exp_m(t\wm) \\
{}             & = & \dte\dt\mid_{t=0}\exp_m(t\wm) \\
{}             & = & \dte d(\exp_m)_0(\wm) \\
{}             & = & \dte\wm
\end{eqnarray*}
Therefore the length of the curve $\wm$ is
\begin{eqnarray*}
l(\wm) & = & \int_0^{2\pi}\parallel\dte\wm\parallel \\
{}     & = & \int_0^{2\pi}\parallel J_{\theta}'(0)\parallel \\
{}     & \leq & 2\pi
\end{eqnarray*}
On the other hand, since $(M, g)$ is a $P_{2\pi}$ manifold at $m$, it
follows that $\wm$ and $\widetilde{\sigma}(-\theta)$ are antipodal points
in $U_mM$ for each $\theta$. Hence $l(\wm)\geq 2\pi$. This shows that
$l(\wm)=2\pi$ and $\wm$ is a great circle. Thus we have proved that
\begin{enumerate}
\item $\parallel J_{\theta}'(0)\parallel=1$ and
$\parallel\jt\parallel=\sin t$. Hence $J_{\theta}$ is an
eigenvectorfield of $\nabla^2u$ with eigenvalue $-u$ for each
$\theta$.
\item $\exp_m : S(0, \pi)\to C$ is a great sphere fibration.
\end{enumerate}
Now we prove that $\eu$ is parallel along $X$ and also an eigensubbundle
of $R(., X)X$ with eigenvalue $1$.

Since $J_{\theta}$ is an eigenvectorfield of $\nabla^2u$
with eigenvalue $-u$ and $\parallel\jt\parallel=\sin t$,
we can write $J_{\theta}(t)=\sin t E(t)$ where $E$ is a
unit eigenvectorfield in $\eu$. Since $E$ is a unit field $E'\perp E$.
But
\begin{eqnarray*}
-J_{\theta}' & = & \nabla_XJ_{\theta} \\
{}           & = & \nabla_{J_{\theta}}X \\
{}           & = & {1\over \minnu}\nabla_{J_{\theta}}\nabu\\
{}           & = & -{u\over\minnu}J_{\theta}
\end{eqnarray*}
i.e., $J_{\theta}'$ is also along $E$.

On the other hand $J_{\theta}'=\cos t E+\sin t E'$. This shows that
$\sin tE'=J_{\theta}'-\cos t E$ is along $E$. Hence $E'=0$ along $X$.
i.e., $E$ is a parallel vectorfield along $X$. This proves that $\eu$ is
parallel along $X$. Hence as in lemma 4 we have that
$\eu$ is also an eigensubbundle of $R(\cdot, X)X$ with
eiegnvalue $1$ and hence $K(E, X)=1$ for $E\in\eu$ of unit norm.

Since $\eu$ is parallel along $X$, $\euu$ is also parallel along $X$. Now
by an easy computation we can show that $R(E, X)X={E\over 4}$ for
$E\in\euu$. Hence $K(E, X)={1\over 4}$ for $E\in\euu$ of unit norm.
This shows that any Jacobi field $J$ describing the variation of a geodesic
$\gamma$ such that $J(0)=0$ and $J'(0)\in\euu$ is given by
$J(t)=2\sin{t\over 2}E(t)$ where $E(t)$ is a parallel field along $X$.
(See also proof of theorem 1 in section 3).

Now we prove the following
\begin{lemma}
$\exp_m : S(0, \pi)\to C$ is congruent to Hopf fibration.
\end{lemma}
{\bf Proof:} The proof is same as that of lemma 5.

\noindent{\bf Proof of theorem 2:} Proof follows along the same lines of
the proof of theorem 3.

\noindent{\bf Concluding Remarks:}
\begin{enumerate}
\item In theorem 2 it is enough to assume that for each $q\in C$ there
exists at least one periodic geodesic $\gamma_q$ of period $2\pi$ from
$m$ to $q$.
\item It appears that an alternative to the hypothesis $P_{2\pi}$ at $m$
could be that diam$(M)=\pi$. Since this forces $m$ to be critical point for
each $d_q$ for $q\in C$.
\item The assumption $P_{2\pi}$ at $m$ is sufficient for our subsequent
work. (See the following remark.)
\item Antonio Ros has proved in \cite{AR} that
if $(M, g)$ is an $n$-dimensional $P_{2\pi}$ manifold  such that the
Ricci tensor $Ric$ and the metric $g$ verify the relation $Ric\geq l g$,
where $l$ is a real constant, then the first eigenvalue $\lambda_1$ of the
Laplacian of $(M, g)$ satisfies the inequality
		$\lambda_1\geq {1\over 3}(2l+n+2)$.
Further the equality
holds iff for any $\lambda_1$ eigenfunction $f$ on $M$ and for any $u$ in
$UM$ we have $f(\gamma_u(t))=A_u\cos t + B_u\sin t +C_u$.
We will discuss in detail the relation between
the results of \cite{AR} and our results elsewhere.
\end{enumerate}

\noindent{\bf Question:} Can one drop the assumption that $(M, g)$ is
$P_{2\pi}$ at $m$? In this context, it may be remarked that it is possible
to construct fibrations of spheres with almost all fibres of diameter
$<2\pi$.

\noindent{\bf Acknowledgements:} We thank Professor K. Grove for
explaining to us the results of \cite{O}, \cite{NS}, which led to the
present formulation of an earlier version of this paper.

\noindent{\em Akhil Ranjan, Department of Mathematics, \\
Indian Institute of Technology,
Bombay- 400 076,  India. \\
e-mail: aranjan@cc.iitb.ernet.in
\vskip.5cm
\noindent G.~Santhanam, School of Mathematics, \\
Tata Institute of Fundamental Research,
Bombay- 400 005,  India. \\
e-mail: santhana@tifrvax.tifr.res.in, santhana@math.tifr.res.in
\vskip.75cm
}
\end{document}